\documentclass[11pt]{amsart}
\usepackage[margin=3cm]{geometry}
\title[First-order convergence for $321$-avoiding permutations]{First-order convergence for $321$-avoiding permutations} 
\author{Alperen \"{O}zdemir} 
\address{Department of Mathematics, KTH Royal Institute of Technology, Stockholm, Sweden} 
\email{alpereno@kth.se} 
\thanks{I would like to thank Michael Damron and Christian Houdré for helpful discussions. 
I also thank the anonymous referee for a careful reading of the manuscript and for numerous insightful comments, which led to significant improvements in the clarity, notation, and rigor of the presentation. 
This work was supported by the Swedish Research Council (VR), grant 2022-03875, and the Knut and Alice Wallenberg Foundation.}
\usepackage{amssymb, amsmath, float, graphicx, qtree, mathtools, pgfplots, tikz, ytableau, xcolor, genyoungtabtikz, mathdots}
\usepackage[colorlinks=true, urlcolor=black,linkcolor=blue, citecolor=blue]{hyperref}
\usepackage{enumerate}

\newtheorem{theorem}{Theorem}[section]
\newtheorem{lemma}{Lemma}[section]

\newtheorem{definition}{Definition}[section]
\newtheorem{example}{Example}[section]

\newtheorem{proposition}{Proposition}[section]

\newcommand{\bbox}{\hfill $\Box$}
\newcommand{\pf}{\noindent {\it Proof:} }

\newcommand{\floor}[1]{\left \lfloor #1 \right \rfloor}
\newcommand{\ceil}[1]{\left \lceil #1 \right \rceil}

\newcommand{\cal}[1]{\mathcal{#1}}
\newcommand{\tn}[1]{\textnormal{#1}}

\usepackage{tikz}
\usetikzlibrary{arrows,shapes}
\usetikzlibrary{trees}
\usetikzlibrary{matrix,arrows}
\usetikzlibrary{positioning}
\usetikzlibrary{calc,through}
\usetikzlibrary{decorations.pathreplacing}
\usepackage{pgffor}

\usetikzlibrary{decorations.pathmorphing}
\usetikzlibrary{decorations.markings}
\tikzset{
    vector/.style={decorate, decoration={snake}, draw},
	provector/.style={decorate, decoration={snake,amplitude=2.5pt}, draw},
	antivector/.style={decorate, decoration={snake,amplitude=-2.5pt}, draw},
    fermion/.style={draw=black, postaction={decorate},
        decoration={markings,mark=at position .55 with {\arrow[draw=black]{>}}}},
    fermionbar/.style={draw=black, postaction={decorate},
        decoration={markings,mark=at position .55 with {\arrow[draw=black]{<}}}},
    fermionnoarrow/.style={draw=black},
    gluon/.style={decorate, draw=black,
        decoration={coil,amplitude=4pt, segment length=5pt}},
    scalar/.style={dashed,draw=black, postaction={decorate},
        decoration={markings,mark=at position .55 with {\arrow[draw=black]{>}}}},
    scalarbar/.style={dashed,draw=black, postaction={decorate},
        decoration={markings,mark=at position .55 with {\arrow[draw=black]{<}}}},
    scalarnoarrow/.style={dashed,draw=black},
    electron/.style={draw=black, postaction={decorate},
        decoration={markings,mark=at position .55 with {\arrow[draw=black]{>}}}},
	bigvector/.style={decorate, decoration={snake,amplitude=4pt}, draw},
}

\begin{document}

\begin{abstract}
We say that a convergence law holds for a sequence of random combinatorial objects if, 
for any first-order sentence $\varphi$, the proportion of objects satisfying $\varphi$ 
converges to a limiting value as the size of the objects tends to infinity. 
In this paper, we show that the convergence law holds for random $321$-avoiding permutations, 
settling an open problem posed in \cite{ABFN22}. Our proof relies on an infinite-dimensional 
version of the Perron-Frobenius theorem.
\end{abstract}

\maketitle

\section{Introduction}

We begin with a class of combinatorial objects and describe it in a language of predicate logic by assigning a model, which is a set equipped with relations and functions, to each object. We assume that there are objects in the class with an arbitrarily large number of elements. We are interested in the limiting fraction of objects that satisfy a given property as the number of elements goes to infinity. More specifically, the properties we are interested in are those that can be expressed in first-order logic. We refer to Section \ref{model} for the precise definitions.

 A question commonly raised in this context is whether the limit described above is always zero or one for a fixed sentence. See the following surveys on logical zero-one law \cite{C89, W93, A18}. A classical result is by Glebskii et al. \cite{GKLT69} and Fagin \cite{F76}. They independently showed that all first-order sentences satisfy zero-one law if we consider the uniform distribution over all models of a given set of elementary relations and functions. The zero-one law appears as a threshold phenomenon for parametric models, such as Erd\H{o}s-Renyi graphs $G(n,p_n)$. 

In first-order logic, one is allowed to quantify over elements only. So its expressive capacity is limited, for instance it does not allow the formulation of sentences involving subsets of the domain as required, for instance for the well-ordering principle or Cantor's theorem. On the other hand, second-order logic allows us to quantify the relations as well as elements, which makes it comparable to set theory \cite{V01}. Unlike first-order logic, it is not complete with respect to standard semantics and it fails the zero-one law \cite{KS85} for uniform distribution over all models of a given language. In fact, the convergence law fails to exist in the general case as one can express the parity (even or odd) of the size of the domain \cite{KV90}.

Within this framework, it is natural to ask whether similar logical limit laws hold for 
specific combinatorial classes. The main problem we consider in this paper is to establish 
a \emph{first-order convergence law} for random $321$-avoiding permutations. Specifically, we ask whether, for any first-order sentence 
$\varphi$, the probability that a uniformly random $321$-avoiding permutation of size $n$ satisfies 
$\varphi$ converges to a limiting value as $n \to \infty$.

There are two different first-order theories to express permutations, TOOB (The Theory of One Bijection) and TOTO (The Theory of Two Orders), which are studied in detail in \cite{ABF20}. The former allows us to express fixed points, cycles, etc., while the latter allows us to articulate the maximum value, adjacent positions, patterns, etc. We will use the latter in this paper. See Example 2.1 for the explicit statement of TOTO. For uniform permutations, it is proved in \cite{FW90} that the limit law does not hold: the authors construct a first-order sentence which distinguishes even and odd permutations. A supplementary demonstration can be found in \cite{F94}.

Let us briefly outline the use of random processes in this context. A finite state space Markov chain is employed for the logical limit laws of random binary words in \cite{L93}. It is recently used in \cite{BK21} to show the convergence law for layered permutations. In a more recent work \cite{MSV23}, the limit laws for permutations under Mallows distributions according to both theories are studied, using countable state space Markov chains and the TOTO representation. Another example is the convergence law for bounded-degree uniform attachment graphs in \cite{MZ22}, where the authors use inhomogeneous Markov chains defined on the subgraphs. In the next section, we show the limitations of Markov chains for $321$-avoiding permutations and describe other processes to address the problem.

To state our main theorem, we first recall the notion of pattern avoidance in permutations.

\begin{definition}[Pattern avoidance] \label{defn:pat}
Let $m \le n$ be positive integers. A permutation $\pi = \pi_1 \cdots \pi_n \in S_n$ 
is said to \emph{contain} a permutation $\tau = \tau_1 \cdots \tau_m \in S_m$ if there exist indices 
$1 \le i_1 < \cdots < i_m \le n$ such that the subsequence $\pi_{i_1}, \ldots, \pi_{i_m}$ 
is in the same relative order as $\tau_1, \ldots, \tau_m$, 
i.e., $\pi_{i_k} < \pi_{i_l}$ if and only if $\tau_k < \tau_l$ for all $1 \le k,l \le m$. We sometimes refer to the shorter permutation $\tau$ as a \emph{pattern}. If $\pi$ does not contain $\tau$, then we say that $\pi$ \emph{avoids} $\tau$, or is \emph{$\tau$-avoiding}.
\end{definition}

We now say that a permutation $\pi$ \emph{satisfies} a first-order sentence $\varphi$ if the property expressed by $\varphi$ holds for $\pi$, denoted by $ \pi \vDash \varphi.$ This notation will be used throughout the paper and is compatible with the more general notion of satisfaction for structures introduced in Section \ref{model}.

We can now state the main result of the paper.
\begin{theorem}\label{mainthm1}
For $\sigma^{321}_n$, a permutation of length $n$ chosen uniformly at random among $321$-avoiding permutations, and for any first-order sentence $\varphi$, the probability that $\sigma^{321}_n$ satisfies $\varphi$ converges to a limit as $n \to \infty$:
\[
\lim_{n \to \infty} \mathbf{P}(\sigma^{321}_n \vDash \varphi) \text{ exists.}
\]
\end{theorem}

The rest of the paper is organized as follows. In Section \ref{sect:prem}, we provide the necessary background on model theory, including definitions of structures, languages, and first-order equivalence. We then introduce the game-theoretical tools and random processes that will be used in the proof of the main theorem. 

In Section \ref{321}, we define a random process on $321$-avoiding permutations and explain how it relates to their logical equivalence classes. Section \ref{CSC} presents the detailed proof of Theorem \ref{mainthm1}, building on the concepts introduced earlier. 

Finally, in Section \ref{sect:ballot}, we derive bounds on certain expressions involving Catalan numbers, which are required in the proof.

\section{Preliminaries} \label{sect:prem}
 
In this section we introduce the model-theoretic concepts and notation needed for the rest of the paper. 
Our goal is to encode permutations as structures in a first-order language so that logical properties 
of permutations can be expressed as first-order sentences. This framework allows us to formulate the 
notion of satisfaction and will later enable us to study convergence laws for random permutations. 
The material presented here provides the background needed for the arguments developed in 
Sections \ref{321} and \ref{CSC}.

\subsection{Model theory}\label{model}
 We refer to \cite{CK90}, \cite{H93} and \cite{M06} for various different presentations of the model theory. Hodges in \cite{H93} humorously conceding to Plato's dialectical reasoning, claims that it is not possible to cover the topic in a linear fashion. Here we provide a concise outline.

A \emph{language} $\cal{L}$ consists of a set of symbols used to describe structures: a set of function symbols $\cal{F}$, a set of relation symbols $\cal{R}$, and a set of constant symbols $\cal{C}$. Each function $f \in \cal{F}$ has a fixed number of arguments, and each relation symbol $R \in \cal{R}$ has a fixed number of arguments. Any of the sets $\cal{F}, \cal{R}, \cal{C}$ may be empty.

The function, relation, and constant symbols of $\cal{L}$ are called \emph{non-logical symbols}, as these are the symbols interpreted by structures. In contrast, the logical symbols, which are present in every language, include negation ($\neg$), equality ($=$), the universal ($\forall$) and existential ($\exists$) quantifiers, and the Boolean connectives ($\vee, \wedge, \Rightarrow, \Leftrightarrow$).

Given a language $\cal{L}$, an \emph{$\cal{L}$-structure} $\cal{M}$ consists of:

\begin{enumerate}
    \item A non-empty set $A$, called the \emph{domain} of $\cal{M}$,
    \item An interpretation $f^\cal{M}: A^{n_f} \to A$ for each function symbol $f \in \cal{F}$,
    \item An interpretation $R^\cal{M} \subseteq A^{n_R}$ for each relation symbol $R \in \cal{R}$,
    \item An element $c^\cal{M} \in A$ for each constant symbol $c \in \cal{C}$.
\end{enumerate}


We introduce \emph{variables} $x_1, x_2, \dots$, which represent elements of a structure. A variable is said to be \emph{free} in a formula if it is not bound by a quantifier, that is, if it does not occur within the scope of $\forall$ or $\exists$. The set of $\cal{L}$-terms is defined inductively as the smallest set containing:
\begin{enumerate}[i)]
    \item every constant symbol of $\cal{L}$,
    \item every variable $x_i$,
    \item expressions $f(t_1,\dots,t_{n_f})$ for each function symbol $f \in \cal{F}$ and terms $t_1,\dots,t_{n_f}$.
\end{enumerate}

An \emph{atomic formula} is either:
\begin{enumerate}[i)]
    \item an equality $t_1 = t_2$ of terms, or
    \item an application $R(t_1,\dots,t_{n_R})$ of a relation symbol $R \in \cal{R}$ to terms $t_1,\dots,t_{n_R}$.
\end{enumerate}
More complex formulas are built from atomic formulas using the logical symbols above.  A \emph{sentence} is a formula with no free variables, and a \emph{theory} is a set of sentences. For a structure $\cal{M}$, we say that $\cal{M}$ is a \emph{model} of a sentence $\varphi$ if $\varphi$ is true in $\cal{M}$, denoted $\cal{M} \vDash \varphi$. If all sentences of a theory $T$ are satisfied by $\cal{M}$, we write $\cal{M} \vDash T$. The \emph{compactness theorem} states that if every finite subset of a first-order theory has a model, then the theory itself has a model. Intuitively, models describe possible worlds, and sentences partition these worlds into true and false cases.

Finally, the \emph{quantifier depth} (or \emph{quantifier rank}) of a formula $\varphi$, 
denoted $\mathrm{qr}(\varphi)$, measures the maximal nesting depth of quantifiers. It is defined inductively as follows:
\begin{enumerate}
    \item $\mathrm{qr}(\varphi) = 0$ if $\varphi$ is atomic,
    \item $\mathrm{qr}(\neg \varphi) = \mathrm{qr}(\varphi)$,
    \item $\mathrm{qr}\!\left(\bigwedge \Phi\right) 
          = \mathrm{qr}\!\left(\bigvee \Phi\right) 
          = \max \{ \mathrm{qr}(\psi) : \psi \in \Phi \}$,
    \item $\mathrm{qr}(\forall x \,\psi) 
          = \mathrm{qr}(\exists x \,\psi) 
          = \mathrm{qr}(\psi) + 1$.
\end{enumerate}

\begin{example} \normalfont 
 (Permutations) We will use TOTO representation where a permutation can be defined by a labelling of elements and their relative positions. Let $S_n$ be the symmetric group on $n$ elements. We write a permutation $\pi \in S_n$ in one-line notation as $\pi=\pi_1\ldots \pi_n.$ Let $[n]:=\{1,2,\ldots,n\}.$ We take $A=[n]$ and define two binary relations, (1) position $(<_P)$ and (2) value $(<_V)$ as follows:
\begin{align*}
i <_P j &\tn{ if and only if } i<j, \\
i <_V j &\tn{ if and only if } \pi_i<\pi_j.
\end{align*}
 These two are linear orders. Therefore a permutation is a structure with $A=[n]$ and $\cal{L}=\{<_P,<_V\}$ two binary relations.
 
  We are interested in the properties of permutations that can be defined as sentences in $\cal{L}$. As an example, consider the sentence
\[\varphi_1=\neg\left[\exists x \exists y \exists z [ (x<_P y)  \, \wedge  \, (y <_P z) \, \wedge \, (y <_V x) \, \wedge \, (z<_V y)]\right],\]
which represents $321$-avoidance. That is to say if $\pi \vDash \varphi_1,$ $\pi$ is a $321$-avoiding permutation.
 A second example is the sentence that ``there exists an inversion", which can be expressed as 
\[\varphi_2 = \exists x \exists y [(x <_P y) \, \wedge \, (y<_V x)].\]
One can also easily show that $\tn{qr}(\varphi_1)=3$ and $\tn{qr}(\varphi_2)=2.$
\end{example}

\subsection{Elementary equivalence}

We will classify the models according to the first-order sentences that they satisfy. The equivalence relations below compare structures defined over the same language.

\begin{definition}\label{eleqclass} \normalfont
Two models $\cal{A}$ and $\cal{B}$ are \emph{elementarily equivalent}, 
denoted $\cal{A} \equiv \cal{B}$, if they agree on all first-order sentences. 
They are \emph{$k$-elementarily equivalent}, denoted $\cal{A} \equiv_k \cal{B}$, 
if they satisfy exactly the same first-order sentences of quantifier depth at most $k$. The equivalence classes induced by $\equiv_k$ will be referred to as the \emph{$k$-equivalence classes}.
\end{definition}

The following result appears as Corollary 3.3.3 in \cite{H93}.
\begin{theorem}
For any two models $\cal{A}$ and $\cal{B},$ $\cal{A} \equiv \cal{B}$ if and only if $\cal{A} \equiv_k \cal{B}$ for all $k \in \mathbb{N}.$
\end{theorem}

Note that if two structures are isomorphic, that is, if there exists a bijection between their domains that preserves the relations and is compatible with the functions, then they satisfy the same first-order sentences. Therefore elementary equivalence is a weaker notion of similarity than isomorphism. For finite structures, the two notions coincide, see Proposition 1.3.19 of \cite{CK90}. An example of non-isomorphic but elementarily equivalent structures is $(\mathbb{R},<)$ and $(\mathbb{Q},<)$. As an example of structures with the same cardinality, consider $\mathbb{Z}$ with the usual order and $\mathbb{Z}^2$ with the lexicographic order, i.e., $(a,n) <_{\tn{lex}}(b,m)$ if and only if $a < b$ or $a=b$ and $n < m$. See also Example 3.3.2 in \cite{Sp01} and Chapter 5 of \cite{H93} for more insight on this topic.

\begin{definition} \normalfont Let two structures $\cal{A}$ and $\cal{B}$ with a common language have domains $A$ and $B$ respectively, and let $S_A \subseteq A$ and $S_B \subseteq B$. A function $g:S_A \to S_B$ is called a \textit{partial isomorphism} if it is a bijection that preserves all relations and functions of $\cal{L}$.

\end{definition}

We will also need information about the number of $k$-equivalence classes.
\begin{theorem}\label{cardlogic}
If $\cal{L}$ contains no function symbols, then there are only finitely many $k$-equivalence classes of $\cal{L}$-structures for each fixed $k \in \mathbb{N}$.
\end{theorem}
\noindent Observe that for any fixed finite set of variables only finitely many atomic formulas can be formed in a relational language. The theorem then follows by a reverse induction argument carried out for binary words; see Lemma 2.3 in \cite{L93}. See also Theorem 2.2.1 in \cite{Sp01} for the corresponding statement in the case of graphs. 

Although the number of $k$-equivalence classes does not depend on the size of the domain of the structures, it can be extremely large. For example, a lower bound is given in \cite{Sp01} in terms of a tower function defined recursively by 
$T(k) = 2^{T(k-1)}$.

\subsection{Ehrenfeucht-Fra\"{i}ss\'{e} Games} \label{sect: ef}

A common method to verify elementary equivalence between two structures is game theoretical. We use a perfect information, sequential and two-player game where the existence of a winning strategy for a player implies the elementary equivalence. We describe it below.

\begin{definition} 
The Ehrenfeucht-Fra\"{i}ss\'{e} game on two sets $A$ and $B$ with $k$ rounds, denoted by $\tn{EF}_k[A,B],$ is played between two players
\begin{align*}
\tn{Player  \textbf{I} } & \, \tn{ a.k.a \, Spoiler, Adam, $\forall$belard}, \\
\tn{Player \textbf{II} }  & \, \tn{ a.k.a  \, Duplicator, Eve. $\exists$loise.}
\end{align*}
At each round, Player \tn{\textbf{I}} first chooses one of the two sets $A$ and $B,$ then chooses an element from that set. Player \tn{\textbf{II}} responds by choosing an element from the other set. Let us denote the element chosen from the set $A,$ by either of the two players, at the end of $l$th stage by $a_l.$ We similarly define $b_l.$  The game is a \textit{win} for Player \tn{\textbf{II}} if there exists a partial isomorphism:
 \[g:\{a_1,\ldots,a_k\}\rightarrow \{b_1,\ldots,b_k\} \tn{ such that }g(a_i)=b_i \tn{ for all } i=1,\ldots,k. \]
\end{definition}
See Theorem 2.4.6 in \cite{M06} for the following: 
\begin{theorem}  \label{gameequiv}
The game $\tn{EF}_k[A,B]$ is a win for Player \tn{\textbf{II}} if and only if $A \equiv_k B.$
\end{theorem}
A worthwhile remark in \cite{Sp01} is that the advantage of the first player is to be able to alternate between two structures.

We will use the following application of the game. First we give a definition. We call a binary relation ``$<$" a \textit{strict linear order} if it is antisymmetric and transitive: 
\[\forall x \forall y  [(x<y) \vee (y<x)] \wedge \neg[(x<y) \wedge (y<x)],\]
\[\forall x \forall y \forall z [(x<y) \wedge (y<z) \Rightarrow (x<z)].\]

\begin{lemma}\tn{(\cite{G84})} \label{gurevich}
Suppose $A$ and $B$ are two structures with a strict linear order, and let $\varphi$ is a sentence of quantifier depth $k.$  $A \equiv_k B$ if and only if $|A|,|B| >2^k$. 
\end{lemma}
\noindent The idea is that, since the order is linear, Player \textbf{II} can always choose an element that preserves the order provided both structures are sufficiently large. Otherwise there will be no interval for Player \textbf{II} to choose from to duplicate at some point before the $k$th stage. 

\subsection{Inhomogeneous processes} \label{nahomojen}

We will consider two growth processes in the rest of the paper. They will first be defined over structures only, then one of them will be extended to the logical equivalence classes as well. They are not time-homogeneous, but we will still be able to derive the convergence to a distribution over the logical equivalence classes. 

This subsection motivates why classical Markov chain techniques are insufficient and prepares the reader for a topological generalization.

Let us start with time-homogeneous Markov chains. If the chain is defined on a finite state space, then the existence of a unique stationary distribution is guaranteed by the Perron-Frobenius theorem under irreducibility and aperiodicity. For countable state spaces, positive recurrence is also required. Countable state space chains have been recently applied in this context \cite{MSV23}. However, in our setting, the transition probabilities cannot be assigned consistently, so the classical Markov chain framework does not apply. Instead, we adopt a topological approach to describe the long-term behavior of the process.

Methods analogous to Markov chains can be found in the dynamical systems literature. An example that we define below is a sequence of distributions given by transition kernels which are not necessarily Markovian. The kernels can be viewed as transfer matrices such as in Section 4.7 of \cite{S11}, but possibly infinite version of them such as in \cite{BJ18}. The idea is to consider a directed graph where the vertices represent the states of a \textit{symbolic chain}, and to look at the frequency of total number of paths terminating at certain vertex in the long-run. A list of references is \cite{V67}, \cite{Sa88} and \cite{GS88}. See Section 7 of \cite{K97} for many other examples and counter-examples. Although we cannot use local arguments, such as Foster's theorem, see Section 2.2 of \cite{FMM95}, we will be able to make inferences on the long-term behavior of the chain.

We outline the definition of a symbolic chain and its properties as given in \cite{K97}. 
We start from a directed graph $\Gamma$ defined on a countable set of vertices $V$ 
together with a binary relation $E \subseteq V \times V$ representing the directed edges. 
We define the adjacency matrix $A$ over $V \times V$ by 
\[
A(i,j)=\mathbf{1}_{(i,j)\in E}.
\]
A symbolic chain is then the sequence of distributions obtained by counting paths 
in $\Gamma$ via powers of the adjacency matrix. 
In particular, the entry $A^n(i,j)$ counts the number of directed paths of length $n$ 
from $i$ to $j$. We call a non-negative matrix \emph{irreducible} if for all states $i,j \in V$ 
there exists $n$ such that $A^n(i,j)>0$. 
We call $d$ a period of the state $i$ if $A^{n+a}(i,i) > 0$ only if 
$n \equiv 0 \pmod d$ for some $a \in \mathbb{N}$. 
The matrix is called \emph{aperiodic} if the common period of all states is $1$.

We next develop tools from the theory of non-negative matrices and directed graphs that will allow us to describe the long-term behavior of the symbolic chains introduced in Section 2.  

The \textit{Perron value} of an irreducible matrix $A$ is
\[\rho(A):=\lim_{n \rightarrow \infty } \sqrt[n]{ A^{n}(i,j)}, \]
which is independent of the choice of $i$ and $j$ by irreducibility. The total number of paths from $i$ to $j$ at $n$th stage is given by 
\begin{equation*}
t_{ij}(n)=A^n(i,j).
\end{equation*}
Let
\[T_{ij}(z)=\sum_{n=0}^{\infty} t_{ij}(n)z^n.\]

Following classical definitions for Markov chains, we classify adjacency matrices based on the growth of path counts and return times. 
We call $A$ \textit{transient} if $T_{ii}(\rho^{-1}(A)) < \infty$ for any $i \in V,$ otherwise it is called \textit{recurrent}. 
To distinguish the chains further, we first define the \textit{first return time} by $f_{ij}(0)=0, f_{ij}(1)=A(i,j)$ and
\[f_{ij}(n+1)=\sum_{k\neq j}A(i,k)f_{ik}(n).\]
Letting $F_i(z)=\sum_{n=1}^{\infty} f_{ii}(n)z^n,$ we call a matrix \textit{positive recurrent} if $F_i'(1)<\infty$ for all $i,$ otherwise it is \textit{null recurrent.}

As in the case of Markov chains, aperiodicity, transience and positive recurrence are class properties; that is, they hold either for all states or for none of them. 

Let us also define the probability simplex for a countable set $V$ as
\begin{equation}\label{simplex}
\Delta^V:=\left\{ (\alpha_i)_{i \in V} \, : \, \sum_{i }\alpha_i =1 \right \},
\end{equation}
which is the boundary of the unit sphere with respect to the $l^1$ norm. We will assume the vertex set is countable. 

The classical Perron-Frobenius theorem guarantees a stationary vector under strong conditions; this will serve as motivation for our topological approach.

\begin{theorem} \label{PF} 
If a non-negative matrix $A$ is irreducible, aperiodic and positive recurrent, then there exists a positive vector $v \in \Delta^V$ such that 
\[v  A=\rho(A)v.\] 
\end{theorem}
\noindent This can be proved by combining theorems 7.1.3 and 7.1.18 in \cite{K97}. 

\vspace{1mm}

In general, the conditions above are not necessary for the existence of a stationary distribution (see example \eqref{expfail} in Section \ref{process}). Let $t(n)$ denote the total number of paths. What is essentially needed for the existence of the stationary distribution of an irreducible chain is 
\[\lim_{n \rightarrow \infty}\frac{t_{ii}(n)}{t(n)} > \delta >0\]
for some $i \in V,$ that is, the ratio of paths leading to some state in the long-run is a positive fraction of the total number of paths. Yet the theorem above requires in addition:
\[\lim_{n \rightarrow \infty}\frac{t(n)}{\rho(A)^n} > C >0.\]

We will give a different formulation to address this issue. The irreducibility corresponds to the connectedness of the graph, and aperiodicity means that the graph is not $k$-partite for any $k\geq 2$ (we call such graphs \textit{non-partite}). To compensate for the positive recurrence, we require the existence of a compact subset of the probability simplex which captures the concentration of the law of the chain on some finite subsets of vertices.  

We now give a topological analogue of Perron-Frobenius suitable for our symbolic chains.

\begin{lemma}\label{TopPF}
Let $\Gamma=(V,E)$ be a locally finite, strongly connected and non-partite directed graph, $\Delta^{V}$ be the probability simplex on the set of vertices and $A$ be the adjacency matrix of $\Gamma.$ Define
\[T:\Delta^{V} \rightarrow \Delta^V \tn{ as } T(w)=\frac{w^T A}{\|w^T A\|_1}.\]
If $T(K)\subseteq K$ for some non-empty, compact and convex $K \subseteq \Delta^{V}$, then there exists a unique $w^* \in K$ such that $\lim_{n \rightarrow \infty} T^n(w_0) =w^*$ for all $w_0 \in \Delta^V.$
\end{lemma}
\noindent The existence follows from Schauder's fixed point theorem \cite[Theorem 2.A]{Z11}, while uniqueness follows from strong connectedness \cite[Theorem 12.2.6]{D07}. \\
We next give a concrete class of examples for the compact set in the lemma, which we will use later. Let $\cal{B}_1$ be a Banach space of sequences with $l^1$ norm, of which $\Delta^V$ is an example. The compact sets of $\cal{B}_1$ are described in Theorem 44.2 of \cite{T67} as the sets which are closed, bounded, and \textit{equismall at infinity}. 

\begin{definition}  
A subset $K$ of $\cal{B}_1$ is \emph{equismall at infinity} if for every $\varepsilon>0,$ there exists $N(\varepsilon)$ such that
 \[\sum_{n>N(\varepsilon)} a_n < \varepsilon \quad \tn{ for all } \, \{a_i\}_{i \in \mathbb{N}} \in K.\]
\end{definition}

Given that $\Delta^X$ is the probability simplex of some countable set $X,$ for any $w\in \Delta^X$ and a real-valued function $f$ on $X$, we define the expected value as
\begin{equation}\label{expectedvalue}
\mathbf{E}_w[f]:= \sum_{x \in X} f(x) w(x).
\end{equation}

Finally, we prove a lemma to identify the compact sets for the cases that we will study. 

\begin{lemma} \label{topexp}
Let $X$ be a countable set with a grading function $f:X \rightarrow \mathbb{N}$ such that $f^{-1}(n)$ is finite for all $n \in \mathbb{N}.$ Then, the set 
\[K=\left \{w \in \Delta^X \, : \, \mathbf{E}_w[f] \leq N \right\}\]
for any $N>0$ is a compact subset of $\Delta^X$.
\end{lemma}

\pf It is clearly bounded and closed with respect to the $l^1$ metric. For every $\varepsilon > 0,$ take $n(\varepsilon)=\ceil{\frac{N}{\varepsilon}}.$ The definition of $K$ implies
\[\sum_{ \{x \in X : f(x)>n_{\varepsilon} \}} w(x)  <\varepsilon.\]
Since $X \setminus \{x \in X : f(x)>n_{\varepsilon} \}$ is finite, we can find a labelling of $X$ to show that $K$ is equismall at infinity.

\bbox

\section{$321$-avoiding permutations} \label{321}

We defined pattern avoiding permutations in Def. \ref{defn:pat}. It is well-known that the numbers of permutations of length $n$ avoiding a pattern of length three are counted by Catalan numbers, 
\begin{equation}\label{catalan}
\cal{C}_n=\frac{1}{n+1}\binom{2n}{n} \sim  \frac{4^n}{\sqrt{\pi} n^{3/2}}.
\end{equation}

See \cite{St15} for a comprehensive survey of Catalan numbers and the many combinatorial structures they enumerate.  

We focus on $321$-avoiding permutations, i.e., permutations $\pi$ of length $n$ such that no three entries, not necessarily adjacent, are in decreasing order: there are no indices $i_1<i_2<i_3$ with $\pi_{i_1}>\pi_{i_2}>\pi_{i_3}.$ Let 
\[\tn{AV}_n(321):=\{\pi \in S_n : \pi \text{ avoids } 321\}, \quad \tn{AV}(321):=\bigcup_{n=1}^{\infty}\tn{AV}_n(321).\]

In \cite{ABFN22}, a first-order limit law is established for $231$-avoiding permutations using their recursive structure. However, the same approach, including Markov chain methods, does not apply to $321$-avoiding permutations. Although $321$- and $231$-avoiding permutations are equinumerous, see Chapter 4 of \cite{K11} for various bijections, they are structurally different in terms of the generating functions of their subclasses \cite{V15}, the distribution of fixed points \cite{E04,MP14,HRS17b}, and their scaling limits \cite{HRS17}. Any other length-3 pattern is equivalent to one of these two up to symmetry.

\begin{figure} 
\begin{tikzpicture} \hspace{0cm}
[xshift=5pt]
\draw (0,0) -- (4,4);
\draw (1,0) -- (5,4);
\fill (0.5,0.5) circle (2pt);
\node at (0.5,0.8) {$1$};
\fill (1.4,1.4) circle (2pt);
\node at (1.4,1.7) {$4$};
\fill (1.7,0.7) circle (2pt);
\node at (1.7,0.4) {$2$};
\fill (2,1) circle (2pt);
\node at (2,0.7) {$3$};
\fill (2.5,2.5) circle (2pt);
\node at (2.5,2.8) {$7$};
\fill (2.8,1.8) circle (2pt);
\node at (2.8,1.5) {$5$};
\fill (3.2,2.2) circle (2pt);
\node at (3.2,1.9) {$6$};
\fill (3.6,2.6) circle (2pt);
\node at (3.6,2.3) {$8$};
\fill (3.9,3.9) circle (2pt);
\node at (3.9,4.2) {$10$};
\fill (4.3,3.3) circle (2pt);
\node at (4.3,3) {$9$};

\node at (2.8,-0.5) {$\pi=1 \, 4 \, 2 \, 3 \, 7 \, 5 \, 6 \, 8 \, 10 \, 9 \in \tn{AV}(321)$};

\fill (7.5,0.2) circle (2pt);
\node at (7.5,0.6) {$1$};
\fill (8,2) circle (2pt);
\node at (8,2.3) {$5$};
\fill (8.5,0.6) circle (2pt);
\node at (8.5,0.9) {$2$};
\fill (9,1.4) circle (2pt);
\node at (9,1.7) {$4$};
\fill (9.5,1) circle (2pt);
\node at (9.5,1.3) {$3$};

\fill (10,4) circle (2pt);
\node at (10,4.3) {$10$};
\fill (10.5,2.4) circle (2pt);
\node at (10.5,2.7) {$6$};
\fill (11.5,3.4) circle (2pt);
\node at (11.5,3.7) {$9$};
\fill (12,3.1) circle (2pt);
\node at (12,3.4) {$8$};
\fill (11,2.8) circle (2pt);
\node at (11,3.1) {$7$};

\node at (9,-0.5) {$\pi=1 \, 5 \, 2 \, 4 \, 3 \, 10 \, 6 \, 7 \, 9 \, 8 \in \tn{AV}(231)$};

 \draw[dashed] (7.3,0)  -- +(0.5,0) -- +(0.5,0.4) -- +(0,0.4) -- cycle;
 \draw[dashed] (8.1,0.5)  -- +(1.5,0) -- +(1.5,1.4) -- +(0,1.4)-- cycle;

 \draw (7.1,-0.1)  -- +(2.6,0) -- +(2.6,2.2) -- +(0,2.2) -- cycle;
 \draw (10.2,2.25)  -- +(2.1,0) -- +(2.1,1.65) -- +(0,1.65) -- cycle;

\draw[dashed] (10.3,2.3)  -- +(1,0) -- +(1,0.6) -- +(0,0.6) -- cycle;
 \draw[dashed] (11.8,3)  -- +(0.4,0) -- +(0.4,0.25) -- +(0,0.25)-- cycle;

\end{tikzpicture}
\caption{Representation of a $321$-avoiding permutation as two increasing sequences, as observed in \cite{M15}. The scaling limit of the two sequences is symmetric, see Theorem 1.2 of \cite{HRS17}. In contrast, $231$-avoiding permutations have a recursive structure, with each box containing another $231$-avoiding permutation.}
\label{slant}
\end{figure}

\subsection{Catalan trees and the insertion process}
To study the limit law for $321$-avoiding permutations, we encode $\tn{AV}(321)$ as a rooted directed infinite graph whose vertices correspond to the permutations themselves. This graph is encoded as a \textit{Catalan tree} (see Figure \ref{fig2}), following \cite{W95}, which captures the recursive insertion structure of the permutations.

\begin{definition}[Catalan trees, \cite{W95}] \label{defn:cat}
Let $\pi \in \tn{AV}_1(321)$ be the root of a rooted directed tree $\mathbf{T}$. Assume that $\pi \in \tn{AV}_n(321)$ labels a vertex at level $n$. For each valid insertion site $i \in \{0,1,\dots,n\}$, define the permutation
\[
\pi^{(i)} = (\pi_1, \dots, \pi_i, n+1, \pi_{i+1}, \dots, \pi_n),
\]
obtained by inserting the new maximum $n+1$ into $\pi$ at position $i+1$, provided the resulting permutation avoids $321$. Each $\pi^{(i)}$ then labels a new child vertex of $\pi$.  

Iterating this construction for all vertices produces a rooted tree $\mathbf{T}$ whose $n$th level consists of all $321$-avoiding permutations of length $n$. This tree is called the \emph{Catalan tree} associated with $321$-avoidance.
\end{definition}

\noindent This construction coincides with the \emph{insertion encoding} of permutations described in \cite{ALR05}, where each permutation is built by inserting a new maximum element into an \emph{active site} in a smaller permutation. Both perspectives are combinatorially equivalent, but the insertion encoding formalism is particularly convenient for algorithmic and probabilistic analysis. In the following, we may omit the adjective “active” and simply write “site” when no ambiguity arises. We will also refer to the inserted new maximum element simply as a “new maximum” once its role is clear.

\begin{definition}[Rank] \label{defn:rank}
Let $G$ be a rooted directed graph. The \emph{rank} of a vertex is the length of the shortest directed path from the root to that vertex.
\end{definition}

\begin{definition}[Truncated Catalan tree] \label{defn:trunc}
For a given $n \ge 1$, let $\mathbf{T}_n$ denote the subgraph of $\mathbf{T}$ containing all vertices of rank at most $n$. We call $\mathbf{T}_n$ the \emph{truncated Catalan tree} at rank $n$.
\end{definition}

\noindent The Catalan tree encodes $\tn{AV}(321)$ combinatorially: each path from the root to a vertex corresponds to the unique sequence of insertions producing a specific permutation. This framework provides a natural way to study recursive constructions, probabilistic growth processes, and bijections with other objects counted by Catalan numbers (see Figure \ref{fig2} for reference).

Next we show that the labelling of the tree defines a matching.
\begin{proposition}
There is a bijection between $\tn{AV}(321)$ and the vertices of $\mathbf{T}$.
\end{proposition}
\pf
Every $\pi \in \tn{AV}_n(321)$ can be traced back to the root $\pi=1$ by recursively removing the largest element, producing only $321$-avoiding intermediate permutations. Conversely, each vertex of $\mathbf{T}$ is uniquely labelled by a permutation in $\tn{AV}(321)$ since the insertion process preserves $321$-avoidance.  
\bbox

\begin{figure}
\centering

\hspace*{-2cm}\Tree[.1 [.$21$ [.$231$ $2341$ $2314$  ]
               [.$213$ $2413$ $2143$ $2134$ ]].$21$
           [.$12$ [.$132$ $1342$ $1324$  ]
               [.$312$ $3412$ $3142$ $3124$ ] 
               [.$123$ $4123$ $1423$ $1243$ $1234$ ]]  .$12$       
]       
       
\vspace*{5mm}

\caption{A Catalan tree representing $321$-avoiding permutations. Each vertex corresponds to a permutation in $\tn{AV}(321)$, and edges indicate insertion of the new maximum into a valid position that preserves $321$-avoidance. The $n$th level contains all permutations of length $n$.}
\label{fig2}
\end{figure}

\subsection{Random processes and ballot numbers}\label{process}

We now connect the combinatorial insertion structure of the Catalan tree with the logical equivalence classes introduced earlier. Each vertex of the tree corresponds to a permutation, and inserting the new maximum produces a directed edge to a permutation of larger size. Thus, the insertion process naturally induces transitions between permutations.

Since permutations are grouped into $k$-elementary equivalence classes, this insertion process also induces transitions between equivalence classes: if $\pi$ and $\pi'$ lie in given classes and $\pi'$ is obtained from $\pi$ by a valid insertion, we obtain a transition between the corresponding classes. Our goal is to analyze these induced transitions and the random processes they generate.

We first study transitions on the Catalan tree itself and identify a key statistic governing the number of active sites. This statistic will later be refined in the proof of the main limit theorem.

The active sites can be described in terms of descents of a permutation.

\begin{definition}\normalfont
Let $\pi$ be a permutation in $S_n$ and define $\pi_0=0.$ $\pi$ is said to have a \emph{descent} at position $i$ for $1 \leq i \leq n-1$ if $ \pi_i > \pi_{i+1}.$ We call a position $i \in \{1,\ldots,n-1\}$ the \emph{rightmost descent} of $\pi$ if there is a descent at position $i$ and there is no other descent at any position $j>i.$ 
\end{definition}

We first describe the active sites.

\begin{lemma} \label{lem:des}
Let $\pi \in \tn{AV}_n(321)$. The active sites for $(n+1)$ are exactly the positions to the right of the rightmost descent of $\pi$ in its one-line notation. If $\pi$ has no descent, then $(n+1)$ may be inserted in any site.
\end{lemma}

This is a standard observation (see, e.g., \cite{W95,ALR05}) and follows directly from the definition of $321$-avoidance.

\begin{lemma} \label{lem:des2}
Let $\pi \in \tn{AV}_{n-1}(321)$, $n\geq 2$, and suppose $\pi$ has $r$ active sites for $n$. Then, after inserting $n$ to obtain $\pi' \in \tn{AV}_n(321)$, the number of sites for $(n+1)$ in $\pi'$ ranges over
\[
\{2,3,\ldots,r,r+1\}.
\]
\end{lemma}
\pf
If $n$ is inserted at the rightmost site, then $\pi'$ has $r+1$ active sites. If $n$ is inserted at the leftmost site, then $\pi'$ has $2$ active sites. Varying the site from right to left produces all intermediate values.
\bbox
      
We are interested in the uniform distribution over $\tn{AV}_n(321)$. The insertion process does not produce this distribution if active sites are chosen uniformly, since different vertices of the Catalan tree have different numbers of branches. The number of branches at a vertex equals one plus the number of elements to the right of the rightmost descent. Let $Q_n$ denote the random variable representing the number of active sites (i.e., one plus the number of elements to the right of the rightmost descent) of a permutation chosen uniformly at random from  $\tn{AV}_n(321)$.

The number of vertices with $r$ branches at level $n$ of the Catalan tree is given by the ballot numbers (see \cite{A01,HP} and OEIS A033184 ballot triangle), namely
\begin{equation}\label{catcoeff}
 q_{n,r} := \frac{r-1}{n} \binom{2n-r}{n-1}, 
 \qquad r=2,\ldots,n+1.
\end{equation}
Their generating function is given in Chapter 1.5 of \cite{FS09}:
\[q_{n,r}=[z^n]\left( \frac{1-\sqrt{1-4z}}{2}\right)^r,\]
from which one obtains the bivariate generating function 
\begin{equation}\label{generatrix}
Q(q,z)=\sum_{r=1}^{\infty}q^r \sum_{n=0}^{\infty} q_{n,r}z^n = \sum_{r=1}^{\infty}q^r \left( \frac{1-\sqrt{1-4z}}{2}\right)^r  = \frac{q(1-\sqrt{1-4z})}{2-q\left(1-\sqrt{1-4z}\right)}.
\end{equation}
The identity
\[\sum_{r=2}^{n+2} q_{n+1,r}=\sum_{r=2}^{n+1} r q_{n,r},\]
which can be observed from the Catalan tree, gives
\begin{equation}\label{expcat}
\mathbf{E}[Q_n]=\frac{\cal{C}_{n+1}}{\cal{C}_n}\rightarrow 4  \tn{ as }n\rightarrow \infty
\end{equation}
by \eqref{catalan}. We can derive the asymptotic distribution of $Q_n$ from the ballot numbers as follows.
\begin{equation} \label{ballotdist}
p(r):=\lim_{n\rightarrow \infty} \mathbf{P}(Q_n=r)=\lim_{n\rightarrow \infty} \frac{q_{n,r}}{\cal{C}_{n}}=\lim_{n\rightarrow \infty}  r\prod_{i=0}^{n-2}\left(1-\frac{r+1}{2n-i} \right) =\frac{r}{2^{r+1}} \ \tn{ for } r=2,3,\ldots
\end{equation}

We will use the tree structure which is generated according to the position of the rightmost descent. We first define two different transition rules on the Catalan tree to address different but related problems. We will also define a third process, which can provide a negative reason for the limitations of Markov chains for these problems.

\begin{enumerate}[(i)]

\item

The first rule is to generate the uniform distribution in a direct way; however, this process does not define a Markov chain. Instead, it corresponds to the symbolic chain described in Section \ref{nahomojen}. Let us take the vertex set to be $V = \mathbb{N} \setminus \{1\}$, which represents the number of children of vertices in the Catalan tree, or equivalently, the set of possible values of the random variable \(Q_n\) (the number of active sites) for all \(n\). Since the number of branches emanating from each vertex depends only on \(Q_n\), this defines a well-posed stochastic process on \(V\).

We define a directed graph on $V$ by putting a directed edge from any vertex $i\in \mathbb{N}$ to the vertices $2,\ldots,i,i+1.$ Suppose the weight of each edge is $1.$ The adjacency matrix of this graph is
\begin{equation} \label{expfail}
(A)_{ij}= \begin{cases}
1 \quad &\tn{ if }j=2,3,\ldots,i+1, \\
0 \quad &\tn {otherwise}.
\end{cases}
\end{equation}

 The Perron value of $A$ is 
\[\rho(A)=\lim_{n\rightarrow \infty} \sqrt[n]{(A^n)_{22}}=\lim_{n\rightarrow \infty} \sqrt[n]{q_{n,2}}=\lim_{n\rightarrow \infty} \sqrt[n]{\cal{C}_{n-1}}=4\]
by the asymptotic formula for the Catalan numbers \eqref{catalan}. This, as expected, agrees with $\mathbf{E}[Q_n].$ We can show that 
\begin{equation} \label{leftev}
l=\left(1,1,\frac{3}{4},\cdots,\frac{n}{2^{n-1}},\cdots \right)
\end{equation}
 is a left eigenvector of $A$. One can also show that the right eigenvector $r=(r_1,r_2,\ldots,)$ satisfies the second-order recurrence relation
\[r_n=4(r_{n-1}-r_{n-2}),\] 
 whose solution gives
\[r=(1,3,8,\ldots,(1+n)2^{n-2},\ldots).\]
It follows that $l \cdot r=\infty,$ which implies that the adjacency matrix is not positive recurrent by Theorem 7.1.3 of \cite{K97}. Even more, one can show that it is transient according to the same lexicon since

\[T_{22}\left(\frac{1}{4}\right)=\sum_{n=0}^{\infty} t_{22}(n) 4^{-n}= \sum_{n=0}^{\infty} \cal{C}_{n-1} 4^{-n} <  \sum_{n=1}^{\infty}n^{-3/2} < \infty\]
 by \eqref{catalan}.

So, although we already identified the distribution of $Q_n$ in \eqref{ballotdist}, its existence does not follow from Lemma \ref{PF}. We will refine this process and define it on a related space to address the logical equivalence in Section \ref{CSC}, then we will use Lemma \ref{TopPF}.

\item The second process that we will use is an inhomogeneous Markov chain that will generate the uniform distribution for a fixed stage. We may use local arguments in this case as in Section \ref{CSC}.

We consider the truncated Catalan tree $\mathbf{T}_n$ and take any vertex $v$ of rank $m \leq n-1$ in $\mathbf{T}_n$ such that there are exactly $r$ directed paths from the vertex to other vertices of rank $(m+1)$ in $\mathbf{T}_n.$ See Def. \ref{defn:rank} and Def. \ref{defn:trunc}. Now let $f(n-m+1,r)$ be the total number of vertices of rank $n$ in $\mathbf{T}_n$ such that there exists a directed path from $v$ to them. It is a simple consequence of Lemma \ref{lem:des} that $f$ does not depend on the choice of the vertex given $n$ (the largest rank in the truncated tree), $m$ (the rank of the initial vertex) and $r$ (the number of branches of the vertex, equivalently the number of admissible insertion positions according to its labelling), so it is well-defined. Observe that $f(2,r)=r$ and $f(n,2)=C_n.$ It can inductively be shown that 
\begin{equation*}
f(n,r)= \binom{2n+r-3}{n-1} - \binom{2n+r-3}{n-2}=\frac{r}{n-1}\binom{2n+r-3}{n-2}.
\end{equation*} 
See also \cite{V14} for various other properties and generalizations of these numbers.

Next we define a finite random process on $\mathbf{T}_n$, whose stationary distribution is the uniform distribution over the vertices of rank $n$. The transition rates will depend on the ranks of vertices and the number of branches of them.

Let $v$ be a vertex of rank $m\leq n-2$ with $r$ branches and $u$ be a vertex connected to $v$ and it is of rank $m+1$ with $i$ branches. It follows from Lemma \ref{lem:des2} that for each possible number of branches $i$, there is a unique vertex of rank $m+1$ connected to $v$, so there is a unique such vertex. To have a simpler expression, we will define probabilities with a decreasing index by taking $l=n-m.$ Then the transition probability from $v$ to $u$ is given by 
\begin{equation*} 
p^l(i|r):=\frac{f(l,i)}{f(l+1,r)}=\frac{i}{r}\frac{l(l+r)\ldots(l+i)}{(2l+r-1)\dots (2l+i-2)} \quad \tn{if } i=2,\ldots,r-1.
\end{equation*}
So we have
\begin{equation} \label{recprob}
p^l(i|r)=\begin{cases}
\frac{r+1}{r}\left(\frac{1}{2}-\frac{r-1}{2(2l+r-1)}\right) \quad &\tn{ if }i=r+1 \\
\frac{1}{4}\left(1-\frac{r-1}{2(2l+r-1)} \right) \left(1+\frac{r+2}{2(2l+r-2)} \right) \quad &\tn{ if } i=r \\
\frac{i}{r \cdot 2^{r-i+2}} \prod_{j=1}^{r+1-i} (1+o\left(\frac{1}{l}\right)) \quad &\tn{ if } i=2,3,\ldots,r-1.
\end{cases}
\end{equation}

 
\item We can extend the process above from truncated trees to  $\mathbf{T}_n$ by letting $l$ go to $\infty.$ The transition rates become uniform and we have indeed a Markov chain with probabilities:
\begin{equation}\label{third}
p(i | r):=\lim_{l \rightarrow \infty} p^l(i|r)=\frac{i}{r \cdot 2^{r-i+2}} \quad \textnormal{ for } i=2,\ldots, r,r+1.
\end{equation}

Consider now the limiting process defined by the transition probabilities in \eqref{third}. Viewed as a stochastic process, 
$\{Q_n\}_{n\ge 1}$ is a random variable taking values in the number of active sites and can be interpreted as a Markov chain on $\mathbb{N} \setminus \{1\}$  with these transition probabilities. However, this Markov chain is not positive recurrent and therefore does not admit a stationary distribution.
 
Formally, if a stationary distribution $\pi$ existed, by definition we would have
\[
\sum_{r=2}^{\infty} \pi(r) \, \mathbb{E}[Q_{n+1}-Q_n \mid Q_n = r] = 0,
\]
because the expected value of $Q_{n+1}$ under $\pi$ must equal the expected value of $Q_n$. But for $r \ge 2$, one can show that
\[
\mathbb{E}[Q_{n+1}-Q_n \mid Q_n = r] = \sum_{i=2}^{r+1} i \, p(i \mid r) - r > 0,
\]
which contradicts this requirement. Hence, no stationary distribution exists.

\end{enumerate}

\subsection{Elementary equivalence and $k$-tails} \label{PEE}

Now we will look at the transitions between logical equivalence classes along with the insertion process. Suppose two permutations in $\tn{AV}(321)$ belong to the same logical equivalence class. For both of them, if we insert the new entry to the rightmost positions according to one-line notation, the new permutations will still have the same equivalence class since the new entries provide no information to make the permutations distinguishable. This can be shown by the EF game, see Section \ref{sect: ef}. Similarly, if it is inserted right next to the rightmost descent, they will have the same equivalence class. However, inserting it into a middle position can possibly give two different $k$-equivalent classes. 

We give an example below that even two $k$-equivalent permutations have the same number of entries right to the rightmost descent, so the same number of insertion positions for the new entries, they are not necessarily $k$-equivalent after inserting the new entry in the same relative positions.

\begin{example}

Suppose $k=3$. Let 
\[
\pi_1 = 3\,\underline{5}\,1\,2\,4\,6\,7\,8 
\quad\text{and}\quad
\pi_2 = 4\,\underline{6}\,1\,2\,3\,5\,7\,8.
\]
The rightmost descents are underlined. 
In each permutation there are six entries to the right of the rightmost descent; 
hence there are seven admissible insertion positions for the new maximum. 
One checks that $\pi_1 \equiv_k \pi_2$ by exhibiting a winning strategy for Player 
\tn{\textbf{II}} in the $k$-round EF game.
 
 It can be shown, however,  that \tn{\textbf{I}} wins the game after the new entries are inserted in the same relative position in both permutations to obtain
\[
\pi_1' = 3\,\underline{5}\,1\,2\,4\,\color{red}9 \, \color{black}6\,7\,8
\quad\text{and}\quad
\pi_2' = 4\,\underline{6}\,1\,2\,3\,\color{red}9 \, \color{black} 5\,7\,8.
\] 
If \tn{\textbf{I}} chooses $\color{red}9$ in $\pi_1'$, then \tn{\textbf{II}} must also choose $\color{red}9$ in $\pi_2'$. Then \tn{\textbf{I}} chooses $5,$ the right-adjacent entry to $\color{red}9$ in $\pi_2'$, and \tn{\textbf{II}} must choose an entry positioned right to $\color{red}9$ in $\pi_1'$. Finally, if \tn{\textbf{I}} chooses the previous rightmost descent, $6,$ in $\pi_2'$, \tn{\textbf{II}} will fail to duplicate. See Figure \ref{anter1} for the picture of a more generic example.
 
\begin{figure} [h!] 
\begin{tikzpicture}
\tikzset{vertex/.style = {shape=circle,draw,minimum size=0.8em}}

 \draw (0,0)  -- +(2,0) -- +(2,2) -- +(0,2) -- cycle;
\draw (2,0) -- +(2,0) -- +(2,2) -- +(0,2);


\fill (2.1,0.1)  circle (1.5pt); \fill (8.1,0.1) circle (1.5pt);
\fill (2.25,0.25)  circle (0.5pt); \fill (8.25,0.25) circle (0.5pt);
\fill (2.35,0.35)  circle (0.5pt); \fill (8.35,0.35) circle (0.5pt);
\fill (2.45,0.45)  circle (0.5pt); \fill (8.45,0.45) circle (0.5pt);

\fill (2.6,0.6)  circle (1.5pt);
 \fill (8.6,0.6) circle (1.5pt);


\fill (0.8,0.7)  circle (1.5pt); \fill (6.8,1) circle (1.5pt);
\draw[dashed] (0.8,0.75) -- (4,0.75); \draw[dashed] (6.8,1) -- (10,1);

\fill (1.3,1)  circle (1.5pt); \fill (7.3,1.25) circle (1.5pt);
\draw[dashed] (1.3,1) -- (4,1); \draw[dashed] (7.3,1.25) -- (10,1.25);

\fill (2.85,0.85)  circle (1.5pt); \fill (8.85,0.85) circle (1.5pt);
\fill (3.15,1.15)  circle (1.5pt); \fill (9.15,1.15) circle (1.5pt);


\fill (3.45,1.5)  circle (0.5pt); \fill (9.45,1.5) circle (0.5pt);
\fill (3.55,1.6)  circle (0.5pt); \fill (9.55,1.6) circle (0.5pt);
\fill (3.65,1.7)  circle (0.5pt); \fill (9.65,1.7) circle (0.5pt);

\fill (3.35,1.4)  circle (1.5pt);
 \fill (9.35,1.4) circle (1.5pt);

\fill (3.85,1.85)  circle (1.5pt); \fill (9.85,1.85) circle (1.5pt);

\fill[red] (3,2.25)  circle (1.5pt); \fill[red] (9,2.25) circle (1.5pt);
\draw[dashed,red] (3,2.25) -- (3,0); \draw[dashed,red] (9,2.25) -- (9,0);

\node at (2,-0.3) {$\pi_A$};
\node at (8,-0.3) {$\pi_B$};

 \draw (6,0)  -- +(2,0) -- +(2,2) -- +(0,2) -- cycle;
 \draw (8,0)  -- +(2,0) -- +(2,2) -- +(0,2);

 \end{tikzpicture}
  \caption{Before the insertion of the new entries (depicted by red dots), Player \tn{\textbf{II}} has a winning strategy if there are sufficiently many entries to the right of the rightmost descents, which are contained in the boxes on the right in the adjacent figures, of both $\pi_A$ and $\pi_B$.}
  \label{anter1}
\end{figure}

\end{example}

The insertion does not respect the equivalence classes in the example above because the new entry can distinguish entries to the right of the rightmost descent that were not initially distinguished; in turn, these could distinguish entries in an antecedent position. To address this issue, we define a decomposition of a certain set of entries, one part of which consists of those located to the right of the rightmost descent. Our intuition is based on the representation of a $321$-avoiding permutation as two increasing sequences, see Figure~\ref{slant}. We thank the referee for pointing out that, when extended to all entries of the permutation, this decomposition is known in the literature as the ``staircase decomposition''; see, for instance, \cite{AABRSW10, AV13, BHPV16}.

We fix a permutation $\pi \in \tn{AV}_n(321)$ for $n\geq 1,$ and define the rightmost descent in a first-order way:
\[\begin{split}
\omega^{\pi}(x):&=  \left[\forall y \, (x>_P y) \,\Rightarrow \, (x>_V y)\right] \, \wedge \, [\exists v \, (v>_P x) \, \wedge \, (x>_V v)]  \\
 & \wedge \, \left[\neg \exists w \left[\forall z \, (w>_P z) \,\Rightarrow \, (w>_V z)\right] \, \wedge \, (w>_P x) \right]
\end{split}\]
Then the entries to the right of the rightmost descent is defined as  
\begin{equation}\label{def:psi1}
\psi_1^{\pi}(x):=  \exists y \, \omega^{\pi}(y)  \Rightarrow (x>_P y).  
\end{equation}
Observe that if the permutation has no descent, then all entries of it satisfy $\psi_1^{\pi}(x).$ See Figure~\ref{anter} for an illustration of the entries satisfying $\psi_1^{\pi}$.

Next we define all entries that have larger value than the smallest entry in $\psi_1^{\pi}(x),$ but does not belong there.
\begin{equation*}
\psi_2^{\pi}(x):= \exists y \, \left( \psi_1^{\pi}(y) \, \wedge \,  x>_V y \right)  \, \wedge \neg \psi_1^{\pi}(x). 
\end{equation*}
For $i \geq 2,$ we first define an auxillary sentence,
\[\varphi_{2i}^{\pi}(x):=  \exists y \,  [(y >_P x) \, \wedge \, (x>_V y)] \, \wedge \left(\neg  \bigvee_{j=1}^{2i-2} \psi_j^{\pi}(x) \right), \]
from which we define other sentences inductively, in order to group other entries of $\pi$ relevant to its equivalence class, as follows:
\begin{equation} \label{defn:psiodd}
\psi_{2i-1}^{\pi}(x):=  \forall y \,  [  \varphi_{2i}^{\pi}(y) \Rightarrow  (x >_P y) ] \, \wedge  \left(\neg  \bigvee_{j=1}^{2i-2} \psi_j^{\pi}(x) \right),
\end{equation}
\begin{equation}  \label{defn:psieven}
\psi_{2i}^{\pi}(x):=  \varphi_{2i}^{\pi}(x) \, \wedge \, \exists y \,  [  \psi_{2i-1}^{\pi}(y) \wedge  (x >_V y)]  \, \wedge \, \neg \psi_{2i-1}^{\pi}(x). 
\end{equation}

Then we define the set of elements of $\pi$ satisfying the sentences above as
\[\psi_i^{\pi}:=\{x \in [n] \, : \, \psi_i(x) \tn{ is true}\}.\]

Observe that $\psi_1^{\bullet}$ is satisfied by the entries in the boxes on the right in Figure \ref{anter1}, while $\psi_2^{\bullet}$ is satisfied by the pair of two entries in the boxes on the left, replacing the superscript ``$\bullet$" by the relevant permutations.



\begin{definition} \tn{($k$-tail of a $321$-avoiding permutation)}
Let $\pi \in \tn{AV}_n(321)$ for some $n\geq 1,$ and let $k\geq 1$ be fixed. We call $\psi^{\pi}=(\psi_1^{\pi},\ldots, \psi_k^{\pi})\subset \{1,2,\ldots,n\}^k$ the \textit{$k$-tail} of $\pi$ and $\psi_i^{\pi}$ the $i$th box of $\psi^{\pi}$ for all $i=1,\ldots,k.$ 
\end{definition}

\begin{proposition}
Let $\pi \in \tn{AV}_n(321)$ for some $n\geq 1,$ and fix $k\geq 1.$ Then the $k$-tail of $\psi^{\pi}=(\psi_1^{\pi},\ldots, \psi_k^{\pi})$ is a decomposition of some non-empty subset of $[n].$
\end{proposition}

\pf For non-emptiness, we show in particular that $\psi_1^{\pi} \neq \emptyset.$ 
If $\pi$ has no descent, then all entries of $\pi$ belong to $\psi_1^{\pi}$, because the premise $\omega^{\pi}(x)$ is false in \eqref{def:psi1}, so that the sentence defining $\psi_1^{\pi}$ is true for all entries. 
On the other hand, if $\pi$ has a descent, there exists a rightmost one, and any entry to the right of it belongs to $\psi_1^{\pi}.$ Therefore, $\psi^{\pi}$ is defined over a non-empty set. 
It is a decomposition, that is, $\psi_i^{\pi} \cap \psi_j^{\pi}=\emptyset$ for $i\neq j.$ 
This follows directly from the definitions, see the last terms in \eqref{defn:psiodd} and \eqref{defn:psieven}.

\bbox

The $k$-tails induce an equivalence relation on $\tn{AV}(321).$ 

\begin{definition}
We call $\pi, \pi' \in \tn{AV}(321)$ \tn{$k$-tail equivalent} if and only if $(\psi_1^{\pi},\ldots, \psi_k^{\pi})=(\psi_1^{\pi'},\ldots, \psi_k^{\pi'}).$ We denote the set of all its equivalence classes by $\Psi.$ 
\end{definition}
When considering the equivalence classes we will drop the reference to $\pi$ and simply write $\psi=(\psi_1,\ldots, \psi_k) \in \Psi.$\\

For our purposes, we label the active sites of a permutation $\pi \in \tn{AV}_n(321)$ from left to right, where each site corresponds to a possible position for inserting the new maximum element. Let $\psi_1^\pi$ be the first box of the $k$-tail, which contains the entries to the right of the rightmost descent. Then the $|\psi_1^\pi|+1$ active sites associated with this set are denoted
\begin{equation} \label{insert}
M_1, M_2, \ldots, M_{|\psi_1^\pi|}, R,
\end{equation}
where $M_i$ represents the new maximum inserted into the $i$th site from the left, and $R$ corresponds to the rightmost site. 

Schematically, we represent sites by empty circles ``$\circ$''. For example, if
\[
\pi = 3\,1\,5\,2\,4\,6,
\]
with rightmost descent underlined as $5>2$, the sites to the right of this descent can be depicted as
\[
3\,1\,5\,\circ\,2\,\circ\,4\,\circ\,6\,\circ{}
\]

There are $4$ sites, and we assign new maximum elements to them as follows:
\begin{itemize}
    \item $M_1$: first $\circ$ (before $2$),
    \item $M_2$: second $\circ$ (between $2$ and $4$),
    \item $M_3$: third $\circ$ (between $4$ and $6$),
    \item $R$: fourth $\circ$ (after $6$, the rightmost site).
\end{itemize}

Formally, let $M$ denote a new maximum inserted into one of the sites, i.e., 
\[M \in \{M_1, M_2, \dots, M_{|\psi_1|}, R\}.\] 
We write the resulting permutation as
\[
\pi' = \pi \,\boxplus\, M,
\]
where $\boxplus$ represents the insertion operation.

Let $\psi$ and $\psi'$ denote the $k$-tails of $\pi$ and $\pi'$ respectively. 
To simplify notation, we also write
\[
\psi \boxplus M := \psi',
\]
so that $(\psi \boxplus M)_i$ refers to the $i$th box of the $k$-tail of the permutation after insertion. 

Then the evolution of the first box of the $k$-tail under this insertion is given by
\[
\psi'_1(x) = (\psi \boxplus M)_1(x) = 
\begin{cases}
\psi_1(x) \vee (x=M), & \text{if } M=R, \\
x>_P M, & \text{if } M \neq R,
\end{cases}
\]
while for $i \ge 2$, $\psi'_i$ is updated according to the inductive $k$-tail construction in \eqref{defn:psiodd}--\eqref{defn:psieven}. In particular, if $M=R$, we have $\psi'_i(x) = \psi_i(x)$ for all $i\ge 2$.  

In this way, the $k$-tail acts as a combinatorial random variable that encodes the local insertion dynamics of the permutation. The evolution of $\psi_1$ coincides with the process on $\mathbb{N} \setminus \{1\}$ describing the number of active sites with stationary distribution \eqref{leftev}, and this connection allows us to track the probabilistic behavior of insertions in a manner consistent with the Catalan tree framework. \\

Now we will show that the evolution of this set under insertion is well-defined with respect to the logical equivalence classes. In particular, we show that if two permutations have the same $k$-tails and the same logical class, then under insertion, the equivalence classes of the new permutations will be the same. Note that the new class can be different from the earlier one. 

\begin{figure} [h!] 
\begin{tikzpicture}
\tikzset{vertex/.style = {shape=circle,draw,minimum size=0.8em}}

 \draw (0,0)  -- +(1,0) -- +(1,1) -- +(0,1) -- cycle;
\draw (0,1) -- +(0,1) -- +(1,1) -- +(1,0);
\draw (1,1) -- +(1,0) -- +(1,1) -- +(0,1);
\draw (1,2) -- +(0,1) -- +(1,1) -- +(1,0);
\draw (2,2) -- +(1,0) -- +(1,1) -- +(0,1);

\fill (0.2,0.2)  circle (1.5pt); \fill[dashed] (5.2,0.2) circle (1.5pt);
\fill (0.6,0.6)  circle (1.5pt); \fill (5.6,0.6) circle (1.5pt);
\fill (0.7,0.7)  circle (1.5pt); \fill (5.7,0.7) circle (1.5pt);

\fill (0.35,1.2)  circle (1.5pt); \fill (5.35,1.2) circle (1.5pt);
\fill (0.5,1.3)  circle (1.5pt); \fill (5.5,1.3) circle (1.5pt);
\fill (0.9,1.75)  circle (1.5pt); \fill (5.9,1.75) circle (1.5pt);

\fill (1.1,1.1)  circle (1.5pt); \fill (6.1,1.1) circle (1.5pt);
\fill (1.3,1.3)  circle (1.5pt); \fill (6.3,1.3) circle (1.5pt);
\fill (1.5,1.45)  circle (1.5pt); \fill (6.5,1.45) circle (1.5pt);
\fill (1.75,1.65)  circle (1.5pt); \fill (6.75,1.65) circle (1.5pt);

\fill (1.6,2.3)  circle (1.5pt); \fill (6.6,2.3) circle (1.5pt);
\fill (1.9,2.6)  circle (1.5pt); \fill (6.9,2.6) circle (1.5pt);

\fill (2.15,2.15)  circle (1.5pt); \fill (7.15,2.15) circle (1.5pt);
\fill (2.5,2.5)  circle (1.5pt); \fill (7.5,2.5) circle (1.5pt);
\fill (2.8,2.8)  circle (1.5pt); \fill (7.8,2.8) circle (1.5pt);

\fill[red] (7.3,3.2)  circle (1.5pt);

 \draw[dashed] (5,0)  -- +(1,0) -- +(1,1) -- +(0,1) -- cycle;
\draw[dashed] (5,1) -- +(0,1) -- +(1,1) -- +(1,0);
\draw[dashed] (6,1) -- +(1,0) -- +(1,1) -- +(0,1);
\draw[dashed] (6,2) -- +(0,1) -- +(1,1) -- +(1,0);
\draw[dashed] (7,2) -- +(1,0) -- +(1,1) -- +(0,1);


 \draw (7.4,2.4)  -- +(1,0) -- +(1,1) -- +(0,1) -- cycle;
 \draw (7.4,2.4)  -- +(-0.75,0) -- +(-0.75,1) -- +(0,1) ;
\draw (7.4,2.4)  -- +(-0.75,0) -- +(-0.75,1) -- +(0,1) ;  
\draw (6.65,2.4)  -- +(0,-0.8) -- +(0.75,-0.8) --+(0.75,0) ;
\draw (6.65,2.4)  -- +(-1.1,0) -- +(-1.1,-0.8) --+(0,-0.8) ;
\draw (5.55,1.6)  -- +(0,-1.1) -- +(1.1,-1.1) --+(1.1,0) ;
\node at (1.8,-0.2) {\Large $\psi$};
\node at (6.9,-0.2) {\Large $\psi \rightarrow \psi'$};
\node at (2.8,2.2) {\small $ \psi_1$};
\node at (1.25,2.8) {\small $ \psi_2$};
\node at (1.75,1.2) {\small $ \psi_3$};
\node at (0.3,1.8) {\small $ \psi_4$};
\node at (0.75,0.2) {\small $ \psi_5$};

 \end{tikzpicture}
  \caption{ Suppose $k=5$. The $k$-tail of $\psi$ evolves into $\psi'$ following the insertion of the new entry. The entries in $\psi_1$ through $\psi_5$ are located in the boxes from the northeast to southwest direction. The red dot represents the new entry. }
  \label{anter}
\end{figure}

\begin{lemma} \tn{(Well-definedness)}
Let $\pi, \sigma \in \tn{AV}(321)$ with a common $k$-tail $\psi \in \Psi$ and $\pi \equiv_k \sigma$ for a fixed $k.$  Then the $k$-equivalence classes of the permutations obtained by insertion depend only on the insertion location, i.e.,
\[\pi \boxplus M \equiv_k \sigma \boxplus M \quad \tn{for} \quad M=M_1,M_2,\ldots, M_{|\psi_1|},R.\]
\end{lemma}

\pf We use the EF game to prove the equivalence. Since the two permutations belong to the same equivalence class, Player \textbf{II} has a winning strategy before inserting the new entry. We describe a winning strategy for Player \textbf{II} after inserting the new entry at the same active site $M \in \{M_1, \dots, M_{|\psi_1|}, R\}$:
\begin{enumerate}
\item If Player \textbf{I} chooses the new entry in either permutation, Player \textbf{II} chooses the new entry in the other.
\item If Player \textbf{I} chooses an entry from the common $k$-tail $\psi$, Player \textbf{II} chooses the same entry in the other permutation.
\item If Player \textbf{I} chooses any entry outside the $k$-tail, Player \textbf{II} responds exactly as in the original winning strategy before the insertion.
\end{enumerate}

We argue that this is indeed a winning strategy. If Player \textbf{I} never chooses the new entry, Player \textbf{II} wins by the original strategy. Otherwise, suppose the new entry is chosen at some stage. If all other moves of Player \textbf{I} belong to $\psi$, Player \textbf{II} can duplicate and win. Otherwise, there exists $1 \le i \le k$ such that no move of Player \textbf{I} belongs to $\psi_i$. In this case, Player \textbf{II} can still duplicate all moves belonging to $\psi_j$ for $j \ge i$. Moves in the remaining boxes $\{\psi_1, \dots, \psi_{i-1}\}$ and $\{\psi_{i+1}, \dots, \psi_k\}$ can be considered independently, since for $s<i<t$, 
\[
(x \in \psi_t) \wedge (y \in \psi_s) \Rightarrow (x>_P y) \wedge (x>_V y),
\]
and Player \textbf{II} has winning strategies for both blocks by the hypothesis.

\subsection{Proof of Theorem \ref{mainthm1}} \label{CSC}

We prove the theorem by encoding the insertion process of $321$-avoiding permutations as a symbolic stochastic process on a state space that records both the logical equivalence class and the $k$-tail configuration. 
We then analyze the directed graph induced by insertions and verify that the hypotheses of Lemma \ref{TopPF} hold on each irreducible component. This yields the existence and uniqueness of a stationary distribution for the process.

Let $\mathcal{K}$ be the set of all elementary equivalence classes, which is finite by Theorem \eqref{cardlogic}, and $\Psi$ be the set of equivalence classes of $k$-tails for a fixed $k$, which is a countable infinite set. We consider the product space $S=\mathcal{K} \times \Psi.$ We will define a random process on this space, which extends the insertion process defined in Section \ref{process}, and apply Lemma \ref{TopPF} to a sequence of vectors in $\Delta^S,$ which denotes the probability simplex of $S$ as in \eqref{simplex}.


Let $\Gamma=(V,E)$ be a directed graph with $V=S$ and $E \subset S \times S,$ the ordered pairs of vertices. We place a directed edge from $(A,\psi)$ to $(B,\varphi)$ if the latter state can be obtained from the former by a single insertion step, which we will denote by
\[(A,\psi) \rightarrow (B,\psi').\]

Let us study the structure of the graph $\Gamma$ to verify the hypotheses of Lemma \ref{TopPF}. We start with a few observations:

\begin{enumerate} [i)]

\item \textbf{Null states}: 
For each tail $\psi \in \Psi$ there exists a subset $\mathcal{K}_{\psi} \subseteq \mathcal{K}$ consisting of logical equivalence classes compatible with $\psi$. 
Vertices $(A,\psi)$ with $A \notin \mathcal{K}_{\psi}$ cannot arise from any $321$-avoiding permutation and therefore correspond to isolated vertices of the graph.\\
 
\item \textbf{Connected components}:
We next show that the directed graph $\Gamma$ is not connected; equivalently, the symbolic chain defined on it is not irreducible.

\begin{proposition}
The chain defined on $S$ by insertion is not irreducible.
\end{proposition}

\pf
Consider the neighbors of the root of the Catalan tree, which are labelled by $12$ and $21$. We claim that no two permutations from these two branches can belong to the same $k$-tail equivalence class for $k \ge 3$.  

Suppose we play the EF game between a permutation $\pi$ from the $12$ branch and a permutation $\pi'$ from the $21$ branch. Let Player \textbf{I} choose the smallest entry in $\pi$ (here, $1$ in $<_V$ order). Then Player \textbf{II} must choose the corresponding smallest entry in $\pi'$. On the next turn, Player \textbf{I} chooses the second smallest entry (here, $2$). In $\pi'$, the relative order of these two entries is inverted ($2 <_P 1$), so Player \textbf{II} cannot select an entry that preserves the partial isomorphism.  

More generally, Player \textbf{II} can only win the $k$-round EF game if the $k$ smallest entries in the two permutations have the same relative order. Otherwise, after the first few moves, a discrepancy in the $<_P$ or $<_V$ order appears, preventing a partial isomorphism. By symmetry, the same argument applies if one considers the leftmost, rightmost, or uppermost $k$ entries.  

Therefore, permutations from different branches of the Catalan tree cannot be $k$-tail equivalent, and the symbolic chain is not irreducible.
\bbox

Nevertheless, we can consider irreducible classes individually and study the process there for the existence of the limit. Observe that for any given $k$-tail $\psi \in \Psi$, one can obtain it from any other tail by inserting sufficiently many new maxima in the positions prescribed by $\psi$. Therefore, the chain is irreducible with respect to the second component of $S$. 

We consider the projection of the chain onto the first coordinate $S$, which runs over a finite set $\cal{K}$. There can be transitionary states, call the set of them $T,$ that is to say logical equivalence classes which do not appear for large permutations. Otherwise, there are closed and irreducible subsets of $S$. So we can partition $\cal{K}$ into those subsets and a transient subset, see Theorem 6.3.4 in \cite {GS20}. So we can write
\[S=T \cup S_1 \cup \cdots \cup S_m\]
where $S_1,\ldots, S_m$ are connected graphs. 

A problem to address is how to assign probabilities to irreducible classes, since ergodicity requires irreducibility. Since the initial state is fixed in our case, say 
$v = (1,0,0,\ldots)$, where the first component corresponds to the logical class and the tail of $\pi=1$, we obtain a unique distribution at every stage. First, we assign zero probability in the limit to classes belonging to $T$.

For each closed irreducible class there exists a vertex of the Catalan tree with the following property: all descendants of that vertex remain in the same irreducible class. 
We call such a vertex the \emph{progenitor} of the class. In other words, once the insertion process reaches this vertex, all subsequent insertions remain within that class.

Let $\cal{A}_i$ denote the set of all permutations descending from the progenitor of the $i$th irreducible class. Recalling that $|\mathrm{AV}_n(321)| = C_n$, we then obtain
\begin{equation} \label{irredstates}
\cal{P}_i:=\lim_{n \rightarrow \infty}\mathbf{P}(\sigma_n^{321} \in S_i)=\sum_{\sigma \in \cal{A}_i} \cal{C}_{|\sigma|}^{-1},
\end{equation}
which denotes the asymptotic probability of the $i$th irreducible class.

Therefore, we can show the existence and the uniqueness of the stationary distribution for a single connected component and take the probability above into consideration at the end. Without loss of generality, we will denote our choice of the connected component by $S_1.$ 

\item \textbf{Non-partiteness:}
Finally, we verify that each irreducible component is aperiodic; equivalently, the graph is not $k$-partite for any $k \ge 2$.

Here we refer to Lemma \ref{gurevich}, the EF game defined on linear structures. 
Note that since we restrict our attention to an irreducible class, all states have the same period. 
See Theorem 6.3.2 of \cite{GS20}, which is stated for Markov chains but still applies here.

Let us take a state $s_0=(A,\psi) \in S_1$ such that $|\psi_1|=2$ and associate it with a permutation $\pi \in \tn{AV}(321)$. 
As a consequence of Lemma \ref{gurevich}, the logical equivalence class will be the same whether we insert $2^{k+1}+1$ or $2^{k+1}+2$ entries at the rightmost site. 
We denote this by
\[
\pi \boxplus (2^{k+1}+1)R \equiv_k \pi \boxplus (2^{k+1}+2)R
\]
according to the notation following \eqref{insert}.

Now we insert a new entry at the $(2^k+1)$st site in the former and the $(2^{k}+2)$nd site in the latter. 
Let
\[
\pi'=\pi \boxplus (2^{k+1}+1)R \boxplus M_{2^k+1} 
\quad \tn{and} \quad 
\pi''=\pi \boxplus (2^{k+1}+2)R \boxplus M_{2^k+2},
\]
and denote their $k$-tails by $\psi'$ and $\psi''$. Observe that 
\[|\psi'_1|= |\psi''_1| > 2^k \, \tn{ and } \,|\psi'_3|,|\psi''_3| > 2^k \]
and for all $y \in \psi_2,$
\[x \in \psi'_3 \vee \psi''_3 \Rightarrow (x>_P y) \wedge (x>_V y). \]  
 So again by Lemma \ref{gurevich} again, we can show that the two permutations $\pi'$ and $\pi''$ have the same logical equivalence class. 
 
From there, we couple two cases in terms of insertion. Eventually, no entry from $\psi'_3$ or $\psi''_3$ will remain and the two tails will be the same. Since the insertions are at the same locations, they will also have the same equivalence classes. Then by irreducibility we can move them both to $s_0$. Therefore,
\[t^{(m)}(s_0),t^{(m+1)}(s_0)>1\]
for some $m>1,$ which implies that there exist paths from $s_0$ to itself of lengths $m$ and $m+1$, that is to say $S_1$ is non-partite. \\

\item \textbf{Compactness:}
At this final step, we verify the compactness condition in Lemma \ref{TopPF}. 
Since $\cal{K}$ is finite, compactness reduces to controlling the second component of $S=\cal{K}\times\Psi$. 
We therefore construct a convex and compact set $K\subset \Delta^V$ consisting of distributions over $\Psi$ and show that it is preserved under sufficiently many applications of $T$.

More precisely, we will prove that for some positive integer $N$,
\[
T^N(K)\subseteq K.
\]
The invariance of $K$ under $T$ will then follow at the end of the proof.

We explicitly state the compact set $K$:
\[
\Pi:=\left\{w \in \Delta^V : \mathbf{P}_w(\psi_1=i)=p_i=\frac{i}{2^{i+1}} 
\ \tn{for } i=2,3,\ldots \right\},
\]
which consists of distributions whose first component $\psi_1$ follows the stationary distribution given in \eqref{ballotdist}.

Then we refer to \eqref{expectedvalue} to define
\begin{equation} \label{kompakt}
K_A:=\left \{ w\in \Pi \, : \,   \mathbf{E}_w(|\psi_1|\cdot |\psi|) \leq A \right \} \subset \Delta^V
\end{equation}
for $A>0$. This set is compact by Lemma \ref{topexp}, since the set of $k$-tails of any fixed length is finite. 
It is also convex by the linearity of $\mathbf{E}_w[\cdot]$ with respect to $w$.

Let $w_0\in K_A$ and set $w_N:=T^N(w_0)$. 
We want to show that $w_N\in K_A$ for suitable constants $A$ and $N$, which will be specified later. Let $\psi^{(0)} \in \Psi$ be sampled from $w_0.$ Conditioning on the initial size of the first box of $w_0,$ we first define the expected length of the $k$-tail after $N$ applications of $T$:
\[C_i:=\mathbf{E}_{w_N}\left[\left|\psi \right| \, : \, \left|\psi_1^{(0)}\right|=i \right].\] 
Then, we define
\[L_i := \mathbf{E}_{w_N}[|\psi| \, : \, |\psi_1|=i]\]
and 
\begin{align*}
E_i: = \mathbf{E}_{w_N}\left[|\psi_1| \, : \, \left|\psi_1^{(0)}\right|=i \right]. 
\end{align*}
In the remainder of the proof we split the expectation
\begin{equation}\label{split}
\mathbf{E}_{w_N}(|\psi|\cdot |\psi_1|) 
= \sum_{i=2}^{s} i p_i L_i + \sum_{i=s+1}^{\infty} i p_i L_i,
\end{equation}
for a suitable cutoff $s\ge2$, and treat the two sums separately. For $i\le s$, we show that if $N$ is large enough, it is very likely that all entries of the original $k$-tail are replaced. 
For $i>s$, the size $|\psi_1|$ decreases approximately linearly along the process.

Concerning the former, we consider the second random process defined in Section \ref{process}. 
Recall that its transitions depend only on $\psi_1$. 
Although the transition probabilities are inhomogeneous, the process is Markovian once the time parameter is fixed.

Let $\psi^{(0)},\psi^{(1)},\ldots$ denote the sequence of $k$-tails generated by this process. We now choose the parameter $l$ in that construction to be the stopping time at which the initial tail has been completely replaced. For a starting tail satisfying $|\psi_1^{(0)}|=r$, define
\[
\tau_r=\min\{t\ge1 : \psi^{(t)}\cap\psi^{(0)}=\emptyset\}.
\]
Thus $\tau_r$ is the first time at which none of the entries of the initial tail remain in the configuration.

We will show that

\begin{lemma} \label{negbin}
For $\varepsilon>0$ and $r=2,3,\ldots$, there exists $N(k,r,\varepsilon)$ such that 
\[
\mathbf{P}(\tau_r \geq N(k,r,\varepsilon)) < \varepsilon.
\]
\end{lemma}

\pf
We distinguish three types of active sites (according to \eqref{insert}): 
(1) the leftmost site, producing the new maximum $M_1$, 
(2) the middle sites, producing $M_i$ for $i=2,\ldots,|\psi_1|$, and 
(3) the rightmost site, producing $R$. 
Accordingly, we refer to the inserted maxima $M_1$, $M_i$ $(i\ge 2)$, and $R$ as the \emph{left}, \emph{middle}, and \emph{right maximum}, respectively.

If the new maximum is inserted at the leftmost site (that is, the inserted point is $M_1$), then the number of insertion sites does not change and the new point belongs to $\psi_2$.  
If the insertion occurs at a middle site (producing $M_i$), then the $k$-tail changes as in Figure~\ref{anter}.  
If the insertion occurs at the rightmost site (producing $R$), then the new point belongs to $\psi_1$ and the remaining components of the tail are unchanged.

Insertion at the rightmost site increases $|\psi_1|$ by $1$, whereas insertion at a middle site reduces $|\psi_1|$ by at least $1$.

\medskip

We now establish a structural property of middle insertions.  
Suppose a new maximum $M_i$ is inserted at a middle site. Assume that there exists a point in the permutation to the left of that site which, at the stage of its insertion, was inserted at the rightmost site; denote that earlier maximum by $R'$. 

We claim that $M_i$ cannot belong to a box that contains any point inserted before $R'$. More precisely, let $<_P$ denote the position order (left-to-right in the permutation) and $<_V$ the value order. Then for any maximum $M$ in the tail $\psi$ such that $P <_V R'$,
\begin{equation} \label{alternate}
\bigl[(M <_P R') \wedge (R' <_P M_i)\bigr]
\;\Rightarrow\;
\bigl[\psi_l(M) \wedge \psi_l(M_i) \Rightarrow \psi_l(R')\bigr]
\quad \text{for all } l=1,\ldots,k.
\end{equation}

To see this, suppose first that $R'$ and $M$ belong to the same box $\psi_l$.  
Since $R'$ was inserted as a rightmost entry at its stage of insertion, it was the largest element at that stage. Hence any point inserted earlier satisfies $M' <_V R'$ and therefore $M' <_P R'$.  

If $M_i <_P M''$ for some point $M''$, then necessarily $R' <_V M''$. By definitions \eqref{defn:psiodd} and \eqref{defn:psieven}, it follows that $M_i$ cannot lie in the same box as $R'$. Consequently it cannot lie in the same box as $M$, since
\[
M <_P R' <_P M_i.
\]

If $R'$ and $M$ belong to different boxes, the same ordering relations yield the claim.

\medskip

Now let us take an initial $k$-tail $\psi^{(0)}$ with $|\psi^{(0)}_1|=r$.  
Once $k$ new boxes are created, the length of $\psi$ agrees with the number of inserted points. Observe that the number of insertions required to create $k$ new boxes depends only on $r$, the number of available insertion sites.  

We now construct a random variable $\tau_r$ that dominates the number of insertions required to replace all entries in the initial $k$-tail $\psi^{(0)}$.

The idea is to force the creation of $k$ new rightmost maxima that are separated from the original tail, and then insert sufficiently many middle maxima to push the initial tail entirely out of $\psi$.

To guarantee this, we consider the following two events:
\begin{enumerate}
\item There occurs at least $k$ times that $M_i$ is immediately followed by $R$ for $i\geq 2.$ This ensures 
\[M_{i} <_{V} R \tn{ and } M_{i} <_{P} R\] 
for $k$ different inserted rightmost maxima. Let the associated stopping time be $\tau_{r,1}.$
\item There placed a middle maximum right to the $k$th rightmost maxima obtained above. So we have a middle maximum $M_i'$ such that
\[R <_{V} M_i' \tn{ and } R <_{P} M_i'\]
for all $R$ created in the first event. Letting $Y$ be the random variable counting the number of insertions of $R$ up until $\tau_{r,1},$ this event is guaranteed to happen by the time $r+Y$ many middle maxima are inserted. Note that we consider here the worst case that the middle maximum is placed at the leftmost possible position every time. Let $\tau_{r,2}$ denote that stopping time that $\tau_{r,1}$ (which is greater than $T$ by definition) middle maxima are inserted. 
\end{enumerate}
So the $k$-tail
\begin{equation} \label{alt}
\psi=\psi^{(0)} \boxplus \cdots \boxplus (M_{i_1} \boxplus R) \boxplus  \cdots  \cdots \boxplus (M_{i_k} \boxplus R) \boxplus \cdots  \boxplus M_{i_T}' \quad \tn{ for }\, i_j\geq 2
\end{equation}
does not have a common element with $\psi^{(0)}.$ Therefore, we have $\tau_r\leq\tau_{r,1} + \tau_{r,2}.$ 
 
Since $M_1$ does not contribute to the number of insertion sites, replacing all occurences of $M_1$ by $R$ will give an upper bound and simplify the computations. So we can specialize on a binary sequence with letters $R$ and $M_i$ for $i\geq 2$. Note that the cases for the transition probabilities in \eqref{recprob} are associated with $R,M_1$ and the middle maxima respectively, so we have
\[\mathbf{P}(M=R) \geq \frac{1}{r+1} \quad \tn{ and } \quad \mathbf{P}(M=M_i \, \tn{ for some } \, i=2,\ldots,r-1 ) \geq \frac{1}{12}.\]

Let $S_i$ be a geometric random variable with parameter $\frac{1}{12}$ for all $i,$ which can be coupled with the waiting time for the insertion of a middle point.

Let $(Y_i)_{i=1}^k$ be a sequence of geometric random variables defined recursively as follows: 
\[
Y_i \sim \text{Geom}\left(\frac{1}{r + \sum_{j=1}^i S_j}\right), \quad i=1,\dots,k,
\] 
where each $Y_i$ counts the number of insertions of rightmost maxima after $S_i$ new middle maxima have been created.  

Then we have
\[
\tau_{r,1} = \sum_{i=1}^k S_i + \sum_{i=1}^k Y_i \equiv_d \sum_{i=1}^k \text{Geom}\left(\frac{1}{12}\right) + \sum_{i=1}^k \text{Geom}\left( \frac{1}{r+\sum_{j=1}^i S_j}\right).
\]

The tail bound of the geometric distribution gives
\[
\mathbf{P}(S_i >  12 c)\leq e^{-c}, \quad i=1,\dots,k,
\]
and from this we obtain
\[
\mathbf{P}\left(Y_i>  rc^i+ 12c^2\frac{c^{i}-1}{c-1}\right) < 2 i e^{-c}, \quad i=1,\dots,k.
\]

Finally, by the union bound over $S_i$ and $Y_j$ for all $i$ and $j$, we have
\[\mathbf{P}(\tau_{r,1} > m_1) < k(k+2)e^{-c}\]
where  $m_1=12kc+ rc^{k+1}+12c^{k+2}.$ 

For the second event, we have 
 \[\tau_{r,2}= \sum_{i=1}^{r+\tau_{r,1}} \tn{Geom}\left(\frac{1}{12} \right),\]
which is a negative binomial distribution. Chernoff's inequality (See Section 2.2 of \cite {BLB03}) gives
\[\mathbf{P}(\tau_{r,2} > N_2) \leq \inf_{t} \mathbf{E}[e^{t \cdot \tau_{r,2}}] e^{-N_2t}. \]
We let $N_1=m_1+r$ and optimize over $t$ to have
 \begin{align*}
 \mathbf{P}(\tau_{r,2} > N_2) \leq \left(1+\frac{1}{N_2-N_1}\right)^{N_2} \left( \frac{N_2-N_1}{12}\right)^{N_1} \left(\frac{11}{12}\right)^{N_2-N_1}.
 \end{align*}
If we take $N_2= 2 \left( \log \left( \frac{1}{\delta} \sqrt{\frac{1}{132}} N_1 \right)  \bigg/ \log \frac{12}{11}\right) N_1,$ then 
\[ \mathbf{P}(\tau_{r,2} > N_2) \leq e^2 \delta^{N_1}. \]
Finally, let $c=\log \left(\frac{2k(k+2)}{\varepsilon}\right)$ and $\delta = \frac{1}{2}.$ Then, for $N=N_1+N_2,$ we have
\[\mathbf{P}(\tau_r > N(k,r,\varepsilon)) < \varepsilon.\]
 This completes the proof of Lemma \ref{negbin}.
 
\bbox

\medskip

Here we complete the compactness argument: we show that, for a suitable choice of $N$ and $A$, the distribution after $N$ applications of $T$ remains in the convex and compact set $K_A \subset \Delta^V$.

Letting $N = \max_{2 \le r \le s} N(k,r,\varepsilon)$, where $N(k,r,\varepsilon)$ is as in Lemma \ref{negbin}, we can bound the sums in \eqref{split} as
\[
\mathbf{E}_{w_N}(|\psi|\cdot |\psi_1|) \le (1-\varepsilon)\sum_{i=2}^{s} p_i \, C_i + \varepsilon \sum_{i=2}^{s} i p_i \, (L_i+N) + \sum_{i=s+1}^{\infty} p_i \, E_i (L_i+N).
\]

Applying \eqref{genball} to the last term and simplifying the others yields
\begin{align*}
\mathbf{E}_{w_N}(|\psi|\cdot |\psi_1|) &\le \max_{2\le i \le s} C_i + \varepsilon (N+A) + \left(\frac{4}{s+1} + \frac{2}{N}\right) A + 4N + 2s \\
&\le \max_{2\le i \le s} C_i + N(4+\varepsilon) + \left(\frac{4}{s+1} + \varepsilon + \frac{2}{N}\right) A + 2s.
\end{align*}

Therefore, for 
\[
A > \left( 1 - \frac{4}{s+1} - \varepsilon - \frac{2}{N}\right)^{-1} \left( N(4+\varepsilon) + 2s + \max_{2\le i \le s} C_i \right),
\]
we have
\[
\mathbf{E}_{w_N}(|\psi|\cdot |\psi_1|) < A.
\]
Choosing $s$ sufficiently large ensures $A>0$, and thus $w_N \in K_A$.

\medskip

\item \textbf{Averaging}: Next, we use an averaging argument to finish the proof of Theorem \ref{mainthm1}. We will show that $T^{N}(w)=w$ implies $T(w)=w$. Let $u_0 = u_N = w$ and define
\[
T(u_i) = u_{i+1}, \qquad i=0,1,\ldots,N_1.
\]
Consider the average
\[
u = \frac{1}{N}\sum_{i=1}^N u_i.
\]
Since $T$ is linear, we have
\[
T(u)
=
\frac{1}{N}\sum_{i=1}^N T(u_i)
=
\frac{1}{N}\sum_{i=1}^N u_{i+1}
=
\frac{1}{N}\sum_{i=1}^N u_i
=
u,
\]
where in the last equality we use $u_N = u_0 = w$. Thus $u$ is a fixed point of $T$, and consequently also of $T^N$.

Since $T$ is irreducible and aperiodic, so is $T^N$. Therefore the fixed point of $T^N$ is unique. Because $T^N(w)=w$, it follows that $u=w$. Hence $w$ is itself a fixed point of $T$.

This completes the proof of Theorem \ref{mainthm1}.

\end{enumerate}

We conclude with two problems:
\begin{enumerate}

\item  The law of the Markov chain in \eqref{third} defines a distribution over $321$-avoiding permutations. Does the convergence law hold in that case? 
\item Are there other applications of Lemma \ref{TopPF} in combinatorial contexts? More specifially, what other combinatorial structures have an inhomogeneous growth but nevertheless have a tractable statistic defined on them? 
\end{enumerate}

\section{Ballot numbers}\label{sect:ballot}
 
In this section, we study the number of vertices of the Catalan tree at a given rank that are descendants of a fixed progenitor vertex. 
Let the progenitor $v$ be of rank $r$ and have $s$ branches. Denote by $q_{n,i}^{(s)}$ the number of vertices connected to $v$ that are of rank $n+r$ and have $i$ branches.  
Observe that for $s=2$ we recover the classical ballot numbers \eqref{catcoeff}.

Considering all possible positions for inserting the largest entry, we obtain the recurrence
\[
q_{n+1,i+1}= \sum_{j=i}^{n+1} q_{n,j},
\]
which can be rewritten as
\[
q_{n+1,i+1}= q_{n+1,i}-q_{n,i-1}.
\]

These ballot numbers can also be expressed as alternating polynomials in Catalan numbers.  
For example:
\begin{align*}
q_{n,2}=q_{n,3}&=\cal{C}_{n-1},\\
q_{n+1,4}&=\cal{C}_{n}-\cal{C}_{n-1},\\
q_{n+1,5}&=\cal{C}_{n}-2\cal{C}_{n-1},\\
q_{n+1,6}&=\cal{C}_{n}-3\cal{C}_{n-1} + \cal{C}_{n-2}, \\
q_{n+1,7}&=\cal{C}_{n}-4 \cal{C}_{n-1} + 3 \cal{C}_{n-2}, \\
q_{n+1,8}&=\cal{C}_{n}-5 \cal{C}_{n-1} +6 \cal{C}_{n-2} - \cal{C}_{n-3}, \\
& \ \vdots
\end{align*}

We define the operator
\[
\Delta_s(a_n):=\sum_{i=0}^{\floor{(s-2)/2}}(-1)^i \binom{s-i-2}{i}a_{n-i},
\]
so that
\begin{equation}\label{deltaop}
q_{n+1,s}=\Delta_s(\cal{C}_n), \quad n\ge s+2.
\end{equation}
Since $\Delta_s$ is linear and its coefficients do not depend on the argument, we also have
\begin{equation} \label{deltacom}
\Delta_s (\Delta_t(a_n))=\Delta_t (\Delta_s(a_n)).
\end{equation}

Next, we introduce the generating function
\[
F_s(q,z)=\sum_{n=1}^{\infty} \sum_{k=2}^{n+s-1} q^{(s)}_{n,k} \ q^kz^n
\]
for the numbers $q^{(s)}_{n,k}$.  

\begin{lemma}\label{coeffeqv}
For $s,k\ge2$ and $n\ge1$, we have
\[
q^{(s)}_{n,k} = \Delta_k(q_{n+s,s}).
\]
\end{lemma}

\pf
Considering all possible insertions given an initial vertex of degree $s+1$, we have
\[
F_s(q,z)=q^sz+z(F_2(q,z)+ \cdots + F_{s+1}(q,z)),
\]
which can be rearranged as
\[
F_{s+1}=\frac{1}{z}(F_s-F_{s-1})+q^s+q^{s+1}.
\]

Iterating this gives
\begin{equation}\label{stotwo}
F_s=\frac{f_s(z)}{z^s}F_2 + g(q,z^{-1}),
\end{equation}
where
\[
\delta_s(z):=\sum_{i=0}^{\floor{(s-2)/2}}(-1)^i \binom{s-i-2}{i}z^i.
\]
The operator $\delta_s(z)$ relates to $\Delta_s$ via
\begin{equation}\label{recip}
\sum_{n=0}^{\infty} \Delta_s(a_n) z^n = \sum_{n=0}^{\infty} a_n \delta_s(z) z^n.
\end{equation}

Since $q^{(2)}_{n,k}=q_{n,k}$, applying \eqref{stotwo} and using \eqref{recip}, \eqref{deltaop}, and \eqref{deltacom} yields
\begin{align*}
F_s(q,z) &= \sum_{n=1}^{\infty} \sum_{k=2}^{n+s-1} q_{n+s,k}\ q^k \delta_s(z) z^n \\
&= \sum_{n=1}^{\infty} \sum_{k=2}^{n+s-1} \Delta_s(q_{n+s,k})\ q^k z^n \\
&= \sum_{n=1}^{\infty} \sum_{k=2}^{n+s-1} \Delta_s(\Delta_k(\cal{C}_{n+s-1}))\ q^k z^n\\
&= \sum_{n=1}^{\infty} \sum_{k=2}^{n+s-1} \Delta_k(\Delta_s(\cal{C}_{n+s-1}))\ q^k z^n \\
&= \sum_{n=1}^{\infty} \sum_{k=2}^{n+s-1} \Delta_k(q_{n+s,s})\ q^k z^n,
\end{align*}
establishing the lemma. \bbox

\medskip
Finally, we consider the expected number of insertion positions after $m$ insertions, starting with an initial value $s$:
\[
\mathbf{E}[X_m^{(s)}] = \frac{[z^m]\partial_q F_s(q,z)|_{q=1}}{[z^m]F_s(1,z)}.
\]
By Lemma \ref{coeffeqv}, this can be expressed as
\[
\mathbf{E}[X_m^{(s)}] = \frac{\sum_{k=2}^{m+s-1} k \Delta_k(q_{m+s,s})}{\sum_{k=2}^{m+s-1} \Delta_k(q_{m+s,s})}.
\]

For $s=2$ we recover the classical Catalan ratio:
\[
\mathbf{E}[X_m^{(2)}] = \mathbf{E}[Q_m] = \frac{\sum_{k=2}^{m+2} k \Delta_k(\cal{C}_m)}{\sum_{k=2}^{m+2} \Delta_k(\cal{C}_m)} = \frac{\cal{C}_{m+1}}{\cal{C}_m}.
\]

In general, using the properties of $\Delta_s$, we have
\[
\mathbf{E}[X_m^{(s)}] = \frac{\Delta_s\left( \sum_{k=2}^{m+s-1} k \Delta_k(\cal{C}_{m+s-1})\right)}{\Delta_s\left( \sum_{k=2}^{m+s-1}  \Delta_k(\cal{C}_{m+s-1})\right)} = \frac{\Delta_s(\cal{C}_{m+s})}{\Delta_s(\cal{C}_{m+s-1})} = \frac{q_{m+s+1,s}}{q_{m+s,s}}.
\]

Using \eqref{catcoeff}, this yields the bound
\begin{equation}\label{genball}
\mathbf{E}[X_m^{(s)}] = \frac{(2m+s+2)(2m+s+1)}{(m+s+1)(m+2)} < 4 + \frac{2s}{m}.
\end{equation}
This estimate will be used in Section \ref{CSC} to control the expected size of $k$-tails under successive insertions.

\bibliographystyle{alpha}
\bibliography{patav}

\end{document}